

\input amstex

\documentstyle{amsppt}

\loadbold

\magnification=\magstep1

\pageheight{9.0truein}
\pagewidth{6.5truein}



\def\annr{\operatorname{ann}_R}
\def\Max{\operatorname{Max}}
\def\Irr{\operatorname{Irr}}
\def\NN{{\Bbb N}}
\def\A{{\bold A}}
\def\GoWa{{\bf 1}}
\def\GoZHa{{\bf 2}}
\def\GoZHb{{\bf 3}}
\def\Let{{\bf 4}}
\def\McRob{{\bf 5}}

\topmatter

\title The Closed-Point Zariski Topology for Irreducible Representations
\endtitle

\rightheadtext{Closed-Point Zariski Topology}

\author K. R. Goodearl and E. S. Letzter \endauthor

\address Department of Mathematics, University of California, Santa Barbara,
CA, 93106 \endaddress

\email goodearl\@math.ucsb.edu \endemail

\address Department of Mathematics, Temple University, Philadelphia,
PA 19122 \endaddress

\email letzter\@math.temple.edu \endemail

\thanks The research of both authors was supported in part by Leverhulme Research
Interchange Grant F/00158/X (UK). It was also partially supported by grants from
the National Science Foundation (USA) and the National Security Agency (USA).
\endthanks


\abstract In previous work, the second author introduced a topology, for
spaces of irreducible representations, that reduces to the classical Zariski
topology over commutative rings but provides a proper refinement in various
noncommutative settings. In this paper, a concise and elementary description
of this refined Zariski topology is presented, under certain hypotheses, for
the space of simple left modules over a ring $R$.  Namely, if $R$ is left
noetherian (or satisfies the ascending chain condition for semiprimitive
ideals), and if $R$ is either a countable dimensional algebra (over a field)
or a ring whose (Gabriel-Rentschler) Krull dimension is a countable ordinal,
then each closed set of the refined Zariski topology is the union of a finite
set with a Zariski closed set.  The approach requires certain auxiliary
results guaranteeing embeddings of factor rings into direct products of simple
modules. Analysis of these embeddings mimics earlier work of the first author
and Zimmermann-Huisgen on products of torsion modules.
 
\endabstract

\endtopmatter

\document

\linespacing{1.2}

\head 1. Introduction \endhead

One of the primary obstacles to directly generalizing commutative algebraic
geometry to noncommutative contexts is the apparent absence of a ``one-sided''
Zariski topology -- that is, a noncommutative Zariski topology sensitive to
left (or right) module theory. In \cite{\Let}, the second author introduced a new
``Zariski like'' topology, on the set $\Irr R$ of isomorphism classes of simple left
modules over a ring $R$ (more generally, on the set of isomorphism classes of
simple objects in a complete abelian category). This topology (see (2.4) for the
precise definition), which we will refer to in this paper as the {\it refined
Zariski topology\/}, satisfies the following properties (see \cite{\Let}): First,
when the ring $R$ is commutative (or PI, or FBN), the topology is naturally
equivalent to the classical Zariski topology on the maximal ideal spectrum of $R$.
Second, each point (i.e., each isomorphism class of simple modules) is closed.
Third, when $R$ is noetherian, the topology is noetherian (i.e., the closed subsets
satisfy the descending chain condition). 

A difficulty in working with the refined Zariski topology is that it is far from
obvious how to identify all the open or closed sets, although a large supply is
guaranteed. In particular, it is easily
seen \cite{\Let, 2.5(i)} that this topology is a refinement of the natural Zariski
topology, in which the closed sets all have the form
$$V(I) \; \colon = \; \{ [N] \in \Irr R \mid I.N=0 \} \; ,$$
where $I$ is an ideal of $R$, and where $[N]$ denotes the isomorphism class of a
simple module $N$. Our aim is to show that, under certain relatively mild
assumptions, the closed sets of the refined Zariski topology are precisely the sets
$$V(I) \; \cup \; F \; ,$$
for ideals $I$ of $R$ and finite subsets $F$ of $\Irr R$. Under these
circumstances, the refined Zariski topology can then be described as the
``point closure'' of the Zariski topology, that is, the coarsest topology
under which the points and the Zariski closed subsets are closed.

The hypotheses under which we obtain the above results consist of a weak noetherian
condition -- specifically, the ACC on semiprimitive two-sided ideals -- and a
countability condition. The latter can be either countable vector space dimension
over a field (satisfied, in particular, by any finitely generated algebra), or
countable Krull dimension (in the sense of Gabriel-Rentschler). The weak noetherian
hypothesis already implies that the Zariski topology on $\Irr R$ is noetherian,
and, in the presence of either countability condition, it follows from an
elementary point set topology argument that the refined Zariski topology is
noetherian.
\smallskip

To briefly describe the contents of this paper: Section 2 presents background and
general results. Section 3 concerns algebras, over a field, of countable vector
space dimension, and Section 4 concerns rings with countable Krull dimension.

Throughout, $R$ denotes an associative ring with identity.

\head 2. Definitions and generalities \endhead

\subhead 2.1 The Zariski topology \endsubhead Let $\Irr R$ denote the
set of isomorphism classes of simple left $R$-modules, and for each $p \in
\Irr R$, let $N_p$ denote a chosen simple left $R$-module in the equivalence
class $p$. For ideals $I$ of $R$, set
$$V(I) \; = \; \{ p \in \Irr R \mid I.N_p = 0 \} \; .$$
(Here and below, {\it ideal\/} means {\it two-sided ideal\/}.) It is not hard to see
that the sets $V(I)$ constitute the closed sets of a topology on $\Irr R$; we will
refer to this topology as the {\it Zariski topology\/} on
$\Irr R$, and we will refer to the subsets $V(I)$ as {\it Zariski closed
subsets\/} of $\Irr R$. Note that each closed set has the form $V(J)$ for some
semiprimitive ideal $J$ of $R$, and that $V(J_1) \supsetneq V(J_2)$ for
semiprimitive ideals $J_1 \subsetneq J_2$. Consequently, the Zariski topology on
$\Irr R$ is noetherian if and only if $R$
satisfies the ACC on semiprimitive ideals. 

It is also not hard to see that when $R$ is
commutative, the Zariski topology on $\Irr R$ is naturally equivalent to the
classical Zariski topology on $\Max R$. On the other hand, when $R$ is a simple
ring, the Zariski topology on $\Irr R$ is trivial (i.e., the only closed sets are
$\varnothing$ and $\Irr R$).

\subhead 2.2 The closed-point Zariski topology \endsubhead Let $Y$ be a set
equipped with a topology $\tau$. The {\it point closure\/} of $\tau$ will
refer to the coarsest topology on $Y$ in which both the points and the
$\tau$-closed subsets are closed. We will refer to the point closure of the
Zariski topology on $\Irr R$ as the {\it closed-point Zariski topology\/}. The
following lemma pins down the point closure of a noetherian topology.

\proclaim{2.3 Lemma} Let $\tau$ be a noetherian topology on a set $Y$, and let
$\tau'$ denote the point closure of $\tau$. Then $\tau'$ is a noetherian topology,
and the $\tau'$-closed subsets of $Y$ are precisely the members of the collection
$$\sigma \; \colon = \; \{ C \cup F \mid \text{$C$ is a $\tau$-closed subset of $Y$,
and
$F$ is a finite subset of $Y$}\} \; .$$
\endproclaim

\demo{Proof} We first show that $\sigma$ is the collection of closed sets for a
topology $\tau''$ on $Y$, from which it will clearly follow that $\tau''= \tau'$.
Obviously $\sigma$ contains $\varnothing$ and $Y$, and $\sigma$ is closed under finite
unions. It remains to show that $\sigma$ is closed under arbitrary intersections. 

Let $\bigl( X_\alpha \bigr)_{\alpha\in A}$ be a nonempty collection of sets
from $\sigma$, say each $X_\alpha= C_\alpha \cup F_\alpha$ where $C_\alpha$ is
$\tau$-closed and $F_\alpha$ is finite. Since $\tau$ is noetherian, the collection
of finite intersections of $C_\alpha$s has a minimum element, that is, there exist
indices $\alpha_1,\dots,\alpha_n\in A$ such that every $C_\alpha$ contains the set
$C := C_{\alpha_1} \cap\cdots\cap C_{\alpha_n}$. Now
$$C\subseteq \bigcap_{\alpha\in A} X_\alpha \subseteq X_{\alpha_1} \cap\cdots\cap
X_{\alpha_n} \subseteq C\cup \bigl( F_{\alpha_1}
\cup\cdots\cup F_{\alpha_n} \bigr) \; ,$$
from which we obtain $\bigcap_{\alpha\in A} X_\alpha= C\cup F$ for some
$F\subseteq F_{\alpha_1} \cup\cdots\cup F_{\alpha_n}$. Since $F$ is necessarily
finite, it follows that $\bigcap_{\alpha\in A} X_\alpha \in \sigma$, as desired.
Thus, $\sigma$ coincides with the set of $\tau'$-closed subsets of $Y$.

To see that $\tau'$ is noetherian, let $X_1\supseteq X_2\supseteq \cdots$ be a
descending chain of $\tau'$-closed subsets of $Y$. We may write each $X_i= C_i\cup
F_i$ where $C_i$ is $\tau$-closed and $F_i$ is finite. Note that since
$X_2\subseteq X_1$, we have $X_2= (C_2\cap C_1)\cup (C_2\cap F_1)\cup F_2$, and so
we may replace $C_2$ and $F_2$ by $C_2\cap C_1$ and $(C_2\cap F_1)\cup F_2$
respectively. Continuing inductively, we see that there is no loss of generality in
assuming that
$C_1\supseteq C_2\supseteq \cdots$. Since $\tau$ is noetherian, there is a positive
integer $n$ such that $C_i=C_n$ for $i\ge n$, and we may delete the first $n-1$
terms of our descending chain. Hence, we may assume that $C_i=C_1$ for all $i$.

Now each $X_i= C_1\sqcup (F_i\setminus C_1)$, and $F_1\setminus C_1\supseteq
F_2\setminus C_2 \supseteq \cdots$. Since $F_1$ is finite, there is a positive
integer $m$ such that $F_i\setminus C_1= F_m\setminus C_1$ for $i\ge m$, and
therefore $X_i= X_m$ for $i\ge m$. \qed\enddemo

\subhead 2.4 The refined Zariski topology \endsubhead Let $\A$ be a complete
abelian category, and assume that the collection $\Irr \A$ of isomorphism
classes of simple objects in $\A$ is a set. (This assumption holds, for instance,
when $\A$ has a generator.) For each $p \in \Irr \A$, let $N_p$ denote a chosen
simple object in the isomorphism class $p$. Following \cite{\Let}, we will say
that a subset $X \subseteq \Irr \A$ is {\it closed\/} provided the isomorphism
class of each simple subquotient of the product
$$\prod_{p \in X}N_p$$
is contained in $X$. In \cite{\Let, 2.3} it is shown that these closed sets are
indeed the closed sets of a topology on $\Irr \A$, and in this paper we will call
the resulting topology the {\it refined Zariski topology}. Note that the points are
closed in this topology. In \cite{\Let} it is proved that this topology is
noetherian when $\A$ has a noetherian generator.

It is not hard to show that the refined Zariski topology on $\Irr R = \Irr {\bold
{Mod}}\text{-}R$ is, indeed, a refinement of the Zariski topology described in
(2.1). In \cite{\Let} it is shown, further, that the refined Zariski topology on
$\Irr R$ coincides with the Zariski topology when $R$ is commutative, or PI, or
FBN; in these cases the refined Zariski topology is again naturally equivalent to
$\Max R$ under the classical Zariski topology.

However, it is also shown in \cite{\Let} that the Zariski and refined Zariski
topologies may differ, even in the case of the first Weyl algebra $A_1(k)$
over a field $k$ of characteristic zero: The refined Zariski topology on $\Irr
A_1(k)$ coincides with the finite complement topology (i.e., the nonempty open sets
are the complements of finite sets), while the Zariski topology on $\Irr A_1(k)$
coincides with the trivial topology. In particular, in this example, the
refined Zariski topology and the closed-point Zariski topology coincide.

\subhead 2.5 Notational assumption \endsubhead Unless otherwise noted, when we
refer to $\Irr R$ as a topological space it is the refined Zariski topology to
which we are referring.

\subhead 2.6 The cofinite product condition \endsubhead We isolate the key
module-theoretic property used in our main results.

Let $I$ be an ideal of $R$. 
We will say that $R$ satisfies the {\it cofinite product condition at
$I$\/} provided the following property holds for all families $\bigl( M_i \bigr)_{i
\in \Omega}$ of pairwise nonisomorphic simple left $R$-modules: {\it If
$$I \; = \; \annr \biggl( \prod_{i \in \Omega}M_i \biggr) \; = \;
\annr \biggl( \prod_{i \in \Omega'}M_i \biggr)$$
for all cofinite subsets $\Omega' \subseteq \Omega$, then there is a left $R$-module
embedding }
$$R/I \; \hookrightarrow \; \prod_{i \in \Omega}M_i \; .$$
(Here, and below, $\annr$ is always used to designate the left annihilator
in $R$.) The condition above always holds when $\Omega$ is finite, given the
convention that the direct product of an empty family of modules is the zero
module. Hence, the condition only needs to be checked for infinite index sets
$\Omega$.

In case the ring $R/I$ is a left Ore domain, the cofinite product condition at $I$
can be rephrased as follows: If $\bigl( M_i \bigr)_{i
\in \Omega}$ is any family of pairwise nonisomorphic simple left $(R/I)$-modules
such that the product $\prod_{i \in \Omega}M_i$ is a torsion module, there exists a
cofinite subset $\Omega' \subseteq \Omega$ such that the partial product $\prod_{i
\in \Omega'}M_i$ is a {\it bounded\/} $(R/I)$-module, i.e., its annihilator in
$R/I$ is nonzero. This is a restricted version (for special families of modules) of
the {\it left productively bounded\/} condition studied by Zimmermann-Huisgen and
the first author in \cite{\GoZHa}. Several of the methods used to establish
productive boundedness in \cite{\GoZHa} can be adapted to prove the cofinite
product condition under analogous hypotheses, as we shall see in Sections 3 and 4.
\smallskip

We now state and prove our main abstract result, which identifies the refined
Zariski topology on $\Irr R$ as the closed-point Zariski topology in the presence of
certain restrictions on $R$.

\proclaim{2.7 Theorem} Let $R$ be a ring with the ascending chain condition on
semiprimitive ideals, and assume that $R$ satisfies the cofinite product condition
of {\rm (2.6)} at all of its semiprimitive prime ideals. Then each closed
subset of $\Irr R$ {\rm(}under the refined Zariski topology{\rm)} is equal to $V(I)
\cup S$, for some semiprimitive ideal $I$ of $R$ and some finite subset $S$ of
$\Irr R$. In particular, the refined Zariski and closed-point Zariski topologies
on $\Irr R$ coincide, and this topology is noetherian.
\endproclaim

\demo{Proof} If $X$ is a closed subset of $\Irr R$ and $I := \bigcap_{p\in X}
\annr(N_p)$, then, since $X$ is also a closed subset of $\Irr (R/I)$, we may work
over the semiprimitive ring $R/I$. Thus, there is no loss of generality in assuming
that $R$ is semiprimitive. Second, by noetherian induction, we may (and will) assume
that the conclusion of the theorem holds for all proper semiprimitive factors of
$R$. To set up a contradiction, we assume that $X$ is a closed subset of
$\Irr R$ that is not equal to the union of some Zariski closed subset with
some finite subset. In particular, $X \ne \Irr R$.

Together, the above assumptions immediately imply that (i) any closed
subset (under the refined Zariski topology) contained in a proper Zariski
closed subset of $\Irr R$ is the union of some Zariski closed subset with some
finite subset; (ii) $X$ cannot be contained in any proper Zariski closed
subset of $\Irr R$; and (iii) any Zariski closed subset corresponding to a
nonzero ideal of $R$ must be a proper subset of $\Irr R$.

To start, suppose there are nonzero ideals $I_1$ and $I_2$ of $R$
for which $I_1I_2 = 0$. Set
$$X_1 \; \colon= \; V(I_1) \cap X \qquad \text{and} \qquad X_2 \;
\colon= \; V(I_2) \cap X \; .$$
Note that $X = X_1 \cup X_2$, that $X_1$ and $X_2$ are (refined Zariski)
closed subsets of $\Irr R$, and that $V(I_1)$ and $V(I_2)$ are proper Zariski
closed subsets of $\Irr R$. By observation (i), $X_1$ and $X_2$ can each be
written as the union of a Zariski closed subset with a finite subset. Hence,
$X$ can be written as a union of a Zariski closed subset with a finite subset,
a contradiction.  It follows that $R$ is prime (and semiprimitive).

Next, set 
$$M \; \colon = \; \prod_{p \in X}N_p \; .$$
Note that $X \subseteq V(\annr M)$, and so $\annr M = 0$ by (ii) and (iii)
above.

Now let $S$ be an arbitrary finite subset of $X$. Set $Y = X \setminus S$
and
$$M' \; \colon = \; \prod_{p \in Y}N_p \; .$$
We claim that the semiprimitive ideal $I := \annr M'$ is zero. Note that $Y
\subseteq V(I)$, and so
$$X \; = \; \bigl( V(I) \cap X \bigr) \cup S \; .$$
Moreover, $V(I) \cap X$ is closed (in the refined Zariski topology), and
is contained in the Zariski closed subset $V(I)$. If $I\ne 0$, then, by (iii) and
(i) above,
$$V(I)\cap X \; = \; V(J) \cup T$$
for some ideal $J$ of $R$ and some finite subset $T$ of $\Irr
R$. But then
$$X = V(J) \cup (T \cup S) \; ,$$
a contradiction. Therefore $\annr M'=0$, as claimed.

By definition (2.4), the isomorphism class of any
simple $R$-module subfactor of $M$ is contained in $X$. However, since $\prod_{p
\in Y}N_p$ is faithful for all cofinite subsets $Y\subseteq X$, the cofinite
product condition (at $0$) implies that $R$ embeds as a left
$R$-module into $M$. Therefore, $X = \Irr R$, another contradiction. Thus, we have
proved that the closed subsets of $\Irr R$, under the refined Zariski topology,
have the desired form.

Since the refined Zariski topology is a topology on $\Irr R$, it follows
immediately that this topology coincides with the point closure of the
Zariski topology. Moreover, since $R$ satisfies the ACC on semiprimitive ideals,
the Zariski topology on $\Irr R$ is noetherian, as noted in (2.1). Therefore, it
now follows from Lemma 2.3 that the refined Zariski topology on $\Irr R$ is
noetherian.
\qed\enddemo

\head 3. Application to algebras of countable vector space dimension \endhead

In this section, we consider algebras of countable vector space dimension over a
field $k$. We show that such algebras satisfy the cofinite product condition at all
ideals, and then we apply Theorem 2.7. In fact, countable dimensional algebras
satisfy a much more general cofinite product condition than that of (2.6), as
follows.

\proclaim{3.1 Proposition} Let $R$ be a $k$-algebra with countable $k$-vector space
dimension, $I$ an ideal of $R$, and $\bigl( M_i \bigr)_{i\in \Omega}$ a family of
left $R$-modules {\rm(}not necessarily simple, and not necessarily pairwise
non-isomorphic{\rm)}. If 
$$I \; = \; \annr \biggl( \prod_{i \in \Omega}M_i \biggr) \; = \;
\annr \biggl( \prod_{i \in \Omega'}M_i \biggr)$$
for all cofinite subsets $\Omega' \subseteq \Omega$, then there is a left $R$-module
embedding 
$$R/I \; \hookrightarrow \; \prod_{i \in \Omega}M_i \; .$$
\endproclaim

\demo{Proof} (Mimics \cite{\GoZHa, Proposition 3.3}.) We may assume, without loss of
generality, that $I=0$. Let $b_1,b_2,\ldots$ be a $k$-basis for $R$, and set 
$$M \; \colon = \; \prod_{i \in \Omega}M_i \; ;$$
we must show that there is an embedding $_RR \hookrightarrow
M$.

For $n=1,2,\ldots$, set $B_n \; := \; kb_1 + \cdots + kb_n$.  The main step in
the proof is to construct a sequence
$$\varnothing = F_0 \subsetneq F_1 \subsetneq F_2 \subsetneq \cdots $$
of finite subsets of $\Omega$ such that, for all $n = 1,2,\ldots$, the
following condition holds: 
\roster
\item"$(*)$" There exists an element $x_n \in \prod_{i \in F_n \setminus
F_{n-1}}M_i$ such that $(\annr x_n ) \cap B_n = 0$.
\endroster

Assume that $F_0,\ldots,F_{n-1}$, for some $n>0$, have been chosen satisfying $(*)$.
Choose a nonzero element $v_1 \in B_n$. Since the module
$$\prod_{i \in \Omega \setminus F_{n-1}}M_i$$
is faithful, we can choose an index $i_1 \in \Omega
\setminus F_{n-1}$ and an element $y_1 \in M_{i_1}$ such that $v_1y_1 \ne 0$. If
$(\annr y_1) \cap B_n = 0$, take $F_n = F_{n-1} \sqcup
\{i_1\}$ and $x_n= y_1$. Otherwise, choose a nonzero element $v_2 \in (\annr
y_1) \cap B_n \subsetneq B_n$, and then choose $i_2 \in \Omega
\setminus (F_{n-1} \sqcup \{i_1\})$ and $y_2 \in M_{i_2}$ such that $v_2y_2 \ne 0$,
using the fact that
$$\prod_{i \in \Omega \setminus (F_{n-1}\sqcup \{ i_1 \})} M_i$$
is faithful. Continue in this manner. Since $B_n$ is finite dimensional, we obtain,
eventually, pairwise distinct indices
$i_1,\ldots,i_t$ in $\Omega\setminus F_{n-1}$ and elements $y_j\in M_{i_j}$ for
$j=1,\ldots,t$ with
$$\big(\annr \{y_1,\ldots, y_t\} \big) \cap B_n \; = \; 0.$$
Set $F_n = F_{n-1} \sqcup \{ i_1, \ldots , i_t\}$, and let
$$x_n \in \prod_{i \in F_n \setminus F_{n-1}} M_i$$
be the unique element such that $(x_n)_{i_j}= y_j$ for $j=1,\ldots,t$.
Then $(\annr x_n) \cap B_n = 0$, and the construction of $\varnothing = F_0
\subsetneq F_1 \subsetneq \cdots$ satisfying $(*)$ follows by induction.

To complete the proof of the proposition, let $x$ be the unique element of $M$ such
that
$$x_i \; = \; \cases (x_n)_i &\quad i\in F_n\setminus F_{n-1} \text{\ for some\ }
n\\ 0 &\quad i\notin \bigcup_{n=1}^\infty F_n \; .  \endcases$$
We see that $(\annr x) \cap B_n = 0$ for all $n = 1,2,\ldots$, and so $\annr x =
0$. Therefore $_RR \cong Rx \subseteq M$. \qed\enddemo

\proclaim{3.2 Theorem} Let $R$ be a $k$-algebra with countable $k$-vector space
dimension, and suppose that $R$ satisfies the ascending chain condition on
semiprimitive ideals. Then each closed subset of $\Irr R$ under the refined Zariski
topology is equal to $V(I) \cup S$, for some semiprimitive ideal $I$ of $R$ and some
finite subset $S$ of $\Irr R$. In particular, the refined Zariski and closed-point
Zariski topologies on $\Irr R$ coincide, and are noetherian.\endproclaim

\demo{Proof} Proposition 3.1 and Theorem 2.7. \qed\enddemo

\head 4. Application to algebras of countable Krull dimension \endhead 

In this section, we consider rings with countable left (Gabriel-Rentschler)
Krull dimension. We again establish the cofinite product condition and apply
Theorem 2.7. As in the previous section, we establish a much more general cofinite
product condition than in (2.6), although here we only establish the condition at
prime ideals. Our approach largely follows
\cite{\GoZHa} and \cite{\GoZHb}. The reader is referred to
\cite{\GoWa} and \cite{\McRob}, for example, for background information on Krull
dimension of modules over noncommutative rings.

\subhead 4.1 Notation \endsubhead Given a left $R$-module $M$, let $L(M)$
denote its lattice of submodules, and let $\kappa(M)$ denote the least ordinal
$\kappa$ such that the dual of the interval $[0,\kappa)$ does not embed in
$L(M)$. When we write $L(R)$ or $\kappa(R)$, we regard $R$ as a left
$R$-module.

\subhead 4.2 \endsubhead Let $M$ be a left $R$-module. It is proved in
\cite{\GoZHb, Theorem} that $\kappa(M)$ is countable if and only if $M$ has
countable Krull dimension (i.e., the Krull dimension of $M$ exists and is a
countable ordinal).

\proclaim{4.3 Lemma} Let $\bigl( M_i \bigr)_{i\in \Omega}$ be a family of left
$R$-modules such that
$$\big| \{i\in \Omega \mid \text{\rm $M_i$ is faithful}\} \big| \; \ge \; \big|
\kappa(R) \big|.$$
Then there is an $R$-module embedding of $_RR$ into $M := \prod_{i\in
\Omega } M_i$.
\endproclaim

\demo{Proof} (Mimics \cite{\GoZHa, Proposition 5.3}.) Assume that there exist
no embeddings of $_RR$ into $M$. Consequently, $\annr y \ne 0$
for all $y\in M$. Set $\kappa= \kappa(R)$, and choose a
subset $\Omega'\subseteq \Omega $ of cardinality $\kappa$ such that $M_i$ is
faithful for all $i\in \Omega'$. Since it suffices to embed $_RR$ into $\prod_{i\in
\Omega'}M_i$, there is no loss of generality in assuming that $|\Omega |=|\kappa|$
and that $M_i$ is faithful for all $i\in \Omega $. Moreover, we may assume that
$\Omega = [0,\kappa)$.

We now inductively construct elements $x_\alpha\in M_\alpha$, for $\alpha\in
\Omega $, and left ideals
$$L_\gamma \; \colon= \; \bigcap_{\alpha\le\gamma} \annr x_\alpha \; ,$$
for $\gamma\in \Omega $, such that $L_\beta > L_\gamma$ for all
$\beta < \gamma < \kappa$. To start, choose an arbitrary $x_0\in M_0$. 

Next, let $0<\gamma<\kappa$, and
assume that $x_\alpha$ has been constructed for all $\alpha<\gamma$. Define
$y\in M$ so that $y_\alpha=x_\alpha$ for $\alpha<\gamma$ while $y_\alpha=0$ for
$\alpha\ge\gamma$. Then,
$$\bigcap_{\alpha<\gamma} \annr x_\alpha \; = \; \annr y \; \ne \; 0 \; .$$
Now choose a nonzero element $r\in \annr y$, and, using the fact that
$M_\gamma$ is faithful, an element $x_\gamma\in M_\gamma$ such that
$rx_\gamma \ne 0$. Observe that
$$L_\gamma \; = \; \bigcap_{\alpha\le\gamma} \annr x_\alpha \; < \;
\bigcap_{\alpha<\gamma} \annr x_\alpha \; \le \; \bigcap_{\alpha\le\beta} \annr
x_\alpha \; = \; L_\beta$$
for all $\beta<\gamma$. This establishes the induction step, and the
construction is complete.

To conclude the proof, observe that the rule $\gamma \mapsto L_\gamma$ defines
a map $[0,\kappa) \rightarrow L(R)$ such that $L_\beta> L_\gamma$ for all
$\beta<\gamma<\kappa$. However, the existence of such a map contradicts the
definition of $\kappa$. The lemma is proved. \qed\enddemo

\proclaim{4.4 Lemma} Let $M$ be a left $R$-module and $\bigl( M_i \bigr)_{i\in
\Omega}$ a family of nonzero submodules of $M$. Set
$$S= \biggl\{ \Omega'\subseteq \Omega \; \biggm| \; \bigcap_{j\in \Omega'} M_j
\ne 0 \biggr\} \,,$$
and assume that
\roster \item"(a)" $\Omega \notin S$, \item"(b)" $S$ is closed under finite
unions, \item"(c)" if $\Omega'$ and $\Omega''$ are any disjoint subsets of
$\Omega $, then either $\Omega'\in S$ or $\Omega''\in S$.
\endroster
Then $\kappa(M)$ is uncountable. \endproclaim

\demo{Proof} This is a restatement of \cite{\GoZHa, Lemma 5.4}.
\qed\enddemo

\proclaim{4.5 Proposition} Let $R$ be a ring with countable left Krull dimension,
$I$ a prime ideal of $R$, and $\bigl( M_i \bigr)_{i\in \Omega}$ a family of left
$R$-modules {\rm(}not necessarily simple, and not necessarily pairwise
non-isomorphic{\rm)}. If 
$$I \; = \; \annr \biggl( \prod_{i \in \Omega}M_i \biggr) \; = \;
\annr \biggl( \prod_{i \in \Omega'}M_i \biggr)$$
for all cofinite subsets $\Omega' \subseteq \Omega$, then there is a left $R$-module
embedding 
$$R/I \; \hookrightarrow \; \prod_{i \in \Omega}M_i.$$
\endproclaim

\demo{Proof} (Mimics \cite{\GoZHa, Theorem 5.5}.) Without loss of generality, we
may assume that $R$ is prime and $I=0$. We must show that $_RR$ embeds into the
module
$$M \; \colon= \; \prod_{i\in \Omega} M_i \; .$$

To start, by (4.2), $\kappa(R)$ is countable. Consequently, if infinitely many
$M_i$ are faithful, then the result follows from Lemma 4.3. Hence,
we may assume there exist at most finitely many $i$ such that $M_i$ is
faithful, and we may remove these indices from $\Omega$ without loss of
generality. In particular, we can assume that $A_i := \annr M_i \ne 0$ for all
$i\in \Omega$. Because $M$ is faithful,
$$\bigcap_{i\in \Omega} A_i=0 \; .$$

Now set 
$$S \; = \; \biggl\{ \Omega' \subseteq \Omega \biggm| \bigcap_{i\in \Omega'} A_i \ne
0 \biggr\} \; .$$
Since $R$ is prime, $S$ is closed under finite unions. Because $\kappa(R)$ is
countable, it then follows from Lemma 4.4 that there exist disjoint subsets
$\Gamma_1,\Lambda_1\subseteq \Omega$ which do not belong to $S$, that is,
$$\bigcap_{i\in \Gamma_1} A_i \; = \; \bigcap_{i\in \Lambda_1} A_i \; = \; 0 \; .$$

Similarly, the collection $\{\Omega' \in S\mid \Omega' \subseteq \Lambda_1\}$
is also closed under finite unions, and so Lemma 4.4 implies that there are
disjoint subsets $\Gamma_2,\Lambda_2\subseteq \Lambda_1$ such that
$$\bigcap_{i\in \Gamma_2} A_i \; = \; \bigcap_{i\in \Lambda_2} A_i \; = \; 0 \; .$$
Continuing inductively, we obtain subsets
$\Gamma_1,\Lambda_1,\Gamma_2,\Lambda_2,\dots \subseteq \Omega$ such that
\roster 
\item"(a)" $\Gamma_n\cap \Lambda_n= \varnothing$ 
\item"(b)" $\Gamma_{n+1} \cup \Lambda_{n+1} \subseteq \Lambda_n$ 
\item"(c)" $\bigcap_{i\in \Gamma_n} A_i= \bigcap_{i\in \Lambda_n} A_i= 0$
\endroster
for all $n$. Note that the sets $\Gamma_1,\Gamma_2,\dots$ are pairwise
disjoint.

By construction, each of the modules 
$$P_n \; \colon = \; \prod_{i\in \Gamma_n} M_i$$
is faithful. Since $|\NN| \ge |\kappa(R)|$, it follows from Lemma 4.3 that there is
a left $R$-module embedding
$$R \; \hookrightarrow \; \prod_{n\in\NN} P_n \; = \; \colon P \; .$$
This completes the proof, because $P$ embeds as a left $R$-module into $M$.
\qed\enddemo

\proclaim{4.6 Theorem} Let $R$ be a ring with countable left Krull dimension, and
suppose that $R$ satisfies the ascending chain condition on semiprimitive ideals.
Then each closed subset of $\Irr R$ under the refined Zariski topology is equal to
$V(I) \cup S$, for some semiprimitive ideal $I$ of $R$ and some finite subset $S$
of $\Irr R$. In particular, the refined Zariski and closed-point Zariski topologies
on $\Irr R$ coincide, and are noetherian. \endproclaim

\demo{Proof} Proposition 4.5 and Theorem 2.7. 
\qed\enddemo

\enddocument